\pdfoutput=1

\RequirePackage{fix-cm}

\documentclass[leqno, fleqn, centertags, 12pt]{article}

\usepackage{latexsym}
\usepackage{amsmath}
\usepackage{amscd}
\usepackage{amssymb}
\usepackage{verbatim}

\usepackage[T1]{fontenc}
\usepackage{latexsym}
\usepackage{manfnt}
\usepackage[sans]{dsfont}
\usepackage{palatino}
\usepackage[sc, osf]{mathpazo} 
\usepackage{eulervm}

\usepackage{fullpage}

\usepackage{graphicx}

\newcommand{\aaa}{\mathbb{A}}
\newcommand{\tto}{\,{\to}\,}

\newcommand{\minus}{\hspace*{0.15em}\mbox{\rule[0.4ex]{0.4em}{0.4pt}}\hspace*{0.15em}}

\def\hss{\hskip.025em\ }
\def\sss{\hskip.05em\ }
\def\dss{\hskip.1em\ }
\def\trs{\hskip.15em\ }
\def\qss{\hskip.02em\ }

\def\hff{{\hskip.025em}}
\def\fff{{\hskip.05em}}
\def\dff{{\hskip.1em}}

\def\qff{{\hskip.2em}}

\def\endss{\hspace*{0.05em}}

\def\ffdot{\hspace*{-0.1em}.\hspace*{0.1em}\ }
\def\dfdot{\hspace*{-0.2em}.\hspace*{0.2em}\ }

\def\dfcom{\hspace*{-0.2em},\hspace*{0.2em}\ }

\newcommand{\nsp}{\hspace*{-0.1em}}
\newcommand{\dnsp}{\hspace*{-0.2em}}

\parskip=\the\bigskipamount

\renewcommand{\thefootnote}{\fnsymbol{footnote}}

\makeatletter
\renewcommand{\@makefntext}[1]{\vspace*{0.5ex}\parindent=0em\noindent
\hspace*{-0.4em}
\hbox to 0.4em{\hss\@makefnmark}\hspace*{0.4em}{#1}
}
\makeatother

\newcounter{mysectionnumber}
\newcounter{myparnum}
\newcounter{mylemmanum}[myparnum]
\newcommand{\mysection}[2]{\setcounter{footnote}{0}
\refstepcounter{mysectionnumber}
\section*{ \textnormal{{\themysectionnumber.} {#1}}}\label{#2}\vspace{\bigskipamount}}  

\newcommand{\mynonumbersection}[2]{
\vspace{-2.0ex}
\section*{{}\hspace*{0.0em}\textnormal{{#1}}}\label{#2}\vspace{\bigskipamount}\vspace{-2.0ex}} 

\newcommand{\mypar}[2]{\refstepcounter{myparnum}
{\vspace{-\bigskipamount} \paragraph{\textit{{\themyparnum. #1}\label{#2}}} \hspace{-0.7em}}}
\renewcommand{\themyparnum}{\arabic{myparnum}}

\newcounter{mysubsubnumber}[myparnum]
\newcommand{\myitpar}[1]{\vspace{-\bigskipamount}\paragraph{\textit{#1}}\hspace{-0.7em}}

\newcommand{\myit}[1]{\textbf{\textit{#1}}\hspace{0.0em}}

\newcommand{\proof}{\vspace{-\bigskipamount}{\paragraph{{\emph{Proof}.\hspace*{0.2em}}}\hspace{-0.7em}}}
\newcommand{\eproof}{ $\blacksquare$}

\begin{document}
\title{{Affine\hspace*{0.1em} planes,\hspace*{0.1em} ternary\hspace*{0.1em} rings,\hspace*{0.1em}
and\hspace*{0.1em} examples\hspace*{0.1em}\\ of\hspace*{0.1em} non-Desarguesian\hspace*{0.1em} planes}}
\date{}
\author{\textnormal{Nikolai\hspace*{0.1em} V.\hspace*{0.1em} Ivanov}}                                                                                                                                                                                                                                                                                        
\renewcommand{\thefootnote}{\fnsymbol{footnote}}

\maketitle

\footnotetext{\hspace*{-0.8em}\copyright\ Nikolai V. Ivanov, 2008, 2014.\trs 
Neither the work reported in this paper, 
nor its preparation were supported by any governmental 
or non-governmental agency, foundation, or institution.}

\vspace*{6ex}

\myit{\hspace*{0em}\large Contents}\vspace*{1ex} \vspace*{\bigskipamount}\\ 
\myit{Preface}\hspace*{0.5em}  \hspace*{0.5em} \vspace*{1ex}\\
\myit{\phantom{1}1.}\hspace*{0.5em} Introduction  \hspace*{0.5em} \vspace*{0.25ex}\\
\myit{\phantom{1}2.}\hspace*{0.5em} Affine planes and ternary rings  \hspace*{0.5em} \vspace*{0.25ex}\\
\myit{\phantom{1}3.}\hspace*{0.5em} Isomorphisms of ternary rings  \hspace*{0.5em} \vspace*{0.25ex}\\
\myit{\phantom{1}4.}\hspace*{0.5em} Isotopisms of ternary rings  \hspace*{0.5em} \vspace*{0.25ex}\\
\myit{\phantom{1}5.}\hspace*{0.5em} Weblen-Wedderburn systems  \hspace*{0.5em} \vspace*{0.25ex}\\
\myit{\phantom{1}6.}\hspace*{0.5em} Near-fields, skew-fields, and isomorphisms \hspace*{0.5em} \vspace*{0.25ex}\\
\myit{\phantom{1}7.}\hspace*{0.5em} Translations  \hspace*{0.5em} \vspace*{0.25ex}\\
\myit{\phantom{1}8.}\hspace*{0.5em} Andr\'{e} quasi-fields  \hspace*{0.5em} \vspace*{0.25ex}\\  
\myit{\phantom{1}9.}\hspace*{0.5em} Conclusion: non-Desarguesian planes  \hspace*{0.5em} \vspace*{1ex}\\
\myit{Bibliographie}\hspace*{0.5em}  \hspace*{0.5em}  \vspace*{0.25ex}

\mynonumbersection{Preface}{preface}

{\small 
The main goal of this paper is to present a detailed self-contained exposition of a part of the theory of affine planes
leading to a construction of affine (or, equivalently, projective) planes not satisfying the Desarques axiom.
Unfortunately, most expositions of the theory of affine and projective planes stop before such a construction.
Perhaps, the reason is that this theory is usually presented as a part of combinatorics, while such constructions 
are in the spirit of the abstract algebra. 

This article is intended to be a complement to the introductory expositions of the theory of affine and projective planes,
and as an easy reading for mathematicians with a taste for abstract algebra (and the geometry of points, lines, and planes).

We start with an axiomatic definition of affine planes 
and show that all of them admit a coordinate system over a so-called \emph{ternary ring}
(in fact, ternary rings are defined in such a way as to make this statement true). 
It is well known that the Desargues axiom for an affine plane is equivalent to the existence of a coordinate system over a skew-field.
Since there are excellent presentations of this equivalence (see, for example, the elegant book \cite{Har} by R. Harstshorne),
we simply define a Desarguesian plane as a plane which admits a coordinate system over a skew-field.
Similarly, since the correspondence between affine and projective planes is very well presented in the literature, 
we do not consider the projective planes in this article.

A novelty of our exposition is the notation $(a,x,b)\longmapsto\langle ax+b \rangle$ for the ternary operation in a ternary ring,
replacing the standard notation $(a,x,b)\longmapsto x\cdot a \circ b$\dfdot
The author hopes that the much more suggestive notation $\langle ax+b \rangle$ will made this beautiful theory accessible to a more wide audience.

}

\vspace*{4\bigskipamount}

\renewcommand{\baselinestretch}{1.01}
\selectfont

\mysection{Introduction}{introduction}

\myitpar{Affine planes.} An\dss \emph{affine plane} $\aaa$ is a set, 
the elements of which are called\dss \emph{points},\sss together with a collection of subsets,
called\dss \emph{lines}, satisfying the following three axioms.
\begin{description}
  \item[A1.] \emph{For every two different points there is a unique line containing them.}
  \item[A2.] \emph{For every line $l$ and a point $P$ not in $l$, there is a unique line containing $P$ and disjoint from $l$.}
  \item[A3.] \emph{There are three points such that no line contains all three of them.}
\end{description}
Two lines are called\dss \emph{parallel}\trs if they are either equal, or disjoint.
Note that being parallel is an equivalence relation. 
Indeed, this relation is obviously reflexive and symmetric.
If two lines $l_1\dff,\dff l_2$ are parallel to a line $l$\dfcom
then the intersection of $l_1$ and $l_2$ is empty by the axiom {\bf A2}, i.e. $l_1$ is parallel to $l_2$\dfdot
 
And \emph{isomorphism} of an affine plane $\aaa$ with an affine plane $\aaa'$ 
is defined as a bijection $\aaa\rightarrow\aaa'$ taking lines to lines. 
Two affine planes $\aaa$\dfcom $\aaa'$ are called\sss \emph{isomorphic},
if there exist an isomorphism $\aaa\rightarrow\aaa'$\dfdot

\myitpar{Affine planes and skew-fields.} A \emph{skew-field}\sss is defined in the same way as a field, except that the commutativity of the multiplication is not assumed.
Skew-fields (in particular, fields) lead to the main examples of affine planes. 
Namely\hff, for a skew-field $K$\dfcom let $\aaa=K^2$ and let $(x,y)$ be the canonical coordinates in $K^2$\dfdot
A \emph{line}\sss in $\aaa$ is defined as subset of $\aaa$ described by an equation having either the form $y=ax+b$ for some $a,b\in K$\dfcom 
or the form $x=c$ for some $c\in K$\dfdot  
An easy exercise shows that $K^2$ with lines defined in this way is indeed an affine plane. 
If an affine plane is isomorphic to $K^2$\dfcom then we say that it is\dss \emph{defined over $K$\dfdot}

The class of affine planes defined over skew-fields can be characterized in purely geometric terms, 
i.e. in terms involving only points and lines. 
Namely, an affine plane $\aaa$ is defined over a skew-field if and only if $\aaa$
satisfies the so-called \emph{Minor \emph{and} Major Desargues axioms}. 
Another way to look at this characterization involves \emph{projective planes} 
(which are not used in this article and by this reason are not even defined).
Every affine plane can be canonically embedded in a projective plane, called its \emph{projective completion}, 
by adding a \emph{line at the infinity}\sss to it. 
In particular, one can construct a projective plane starting from a skew-field $K$\dfdot 
A projective plane constructed in this way is said to be \emph{defined over $K$\dfdot}
A projective plane is defined over a skew-field 
if and only if it satisfies the \emph{Desargues axiom} for projective planes. 
Also, the projective completion of an affine plane $\aaa$
is defined over a skew-field $K$\dfcom
i.e. can be obtained from some affine plane by adding the line at infinity
(which may be different from the original line at infinity added to $\aaa$\dnsp)  
if and only $\aaa$ is defined over the same skew-field $K$\dfdot

In the present paper these characterizations serve only as a justification of the following term. 
Namely, an affine plane is said to be \emph{non-Desarguesian} if it is not defined over a skew-field. 
Our main goal is to construct examples of non-Desarguesian planes.

\myitpar{Prerequisites.} The prerequisites for reading this article are rather modest. 
It is not even strictly necessary to be familiar beforehand with the notion of an affine plane. 
But the reader is expected to be familiar with the notions of rings and fields. 
At some places we speak about vector spaces over skew-fields; 
without much loss the reader may assume that these skew-fields are actually fields.
In Section \ref{andre} we use one basic result from the Galois theory, but it can be well taken on faith. 
Mainly, only a taste for abstract algebra is expected, especially in Section \ref{andre}.

\myitpar{The organization of the paper.}  In Section \ref{ternary} we introduce the notion of coordinates in an affine plane,
not necessarily defined over a skew-field. 
The coordinates of a point are taken from a set with a ternary operation, called a \emph{ternary ring}. 
Conversely, every ternary ring defines an affine plane, as explained in Section \ref{ternary}. 
A source of difficulties is the fact that isomorphic affine planes can be coordinatized by non-isomorphic ternary rings. 
They are isomorphic under an obvious additional condition; this is discussed in Section \ref{isomorphisms}. 
In Section \ref{isotopisms} we discuss a weaker notion of an isomorphism for ternary rings 
(namely, the notion of an {\em isotopism}), which is better related to isomorphisms of affine planes. 
But, in fact, this notion is not needed for our main goal, namely, for construction of non-Desarguesian planes,
and Section \ref{isotopisms} may be skipped without loss of the continuity. 
In Section \ref{quasi-fields} we introduce the most tractable class of ternary rings, namely, the class of
\emph{Veblen-Wedderburn systems}, also called \emph{quasi-fields}. 
Like the fields, they are sets with two binary operations, but satisfying only a fairly weak version of the axioms of a field. 
Sections \ref{near-fields} and \ref{translations} are devoted to two different ways to prove that an affine plane is not defined over a skew-field. 
Section \ref{andre} is devoted to a construction of quasi-fields not isomorphic to a skew-field. 
Finally, in Section \ref{conclusion} we combine the results of the previous sections in order to construct non-Desarguesian planes. 
The main ideas are contained in Sections \ref{ternary}, \ref{quasi-fields}, and \ref{andre}.

\myitpar{Further reading.} There are several excellent books exploring deeper these topics, 
in particular, exploring the role of the Desargues axioms. 
For a systematic introduction to the theory of affine and projective planes the reader may turn 
to the classical unsurpassed books by E. Artin \cite{Ar} (see \cite{Ar}, Chapter II),
M. Hall \cite{Ha1} (see \cite{Ha1}, Chapter 20), and R. Hartshorne \cite{Har}. 
The book by Hartshorne is the most elementary one of them. 
The reader is not assumed even to be familiar with the notions of a group and of a field; 
in fact, the affine and projective geometries are used to motivate these notions. 
The book by M. Hall (or, rather, its last chapter, which actually does not depend much on the previous ones) 
gives an in depth exposition directed to mature mathematicians of the theory of projective planes and its connections with algebra. 
E. Artin's elegant exposition is written on an intermediate level between books of R. Hartshorne and M. Hall.
All these books present in details the characterization of planes defined over skew-fields  in terms of Desargues axioms. 

A lot of books in combinatorics discuss some elementary parts of the theory of affine and projective planes, 
but very rarely include a construction of a non-Desarguesian plane. 
An exception is another M. Hall's classics, namely \cite{Ha2}. 
We followed \cite{Ha2} in that we deal with affine planes and not with projective ones, 
and in the way we coordinatize affine planes in Section \ref{ternary}.

The book of D. Hughes and F. Piper \cite{HP} is, probably, 
the most comprehensive exposition of the theory of projective planes 
(the study of which is essentially equivalent to the study of affine planes). 
The state of the art as of 2007 is discussed by Ch. Weibel \cite{W}.

\mysection{Affine planes and ternary rings}{ternary}

Let $\aaa$ be an affine plane. We start by introducing some sort of cartesian coordinates in $\aaa$\dfdot
We follow the approach of M. Hall (see \cite{Ha2}, Section 12.4). 
Then we use these coordinates in order to define a ternary operation on any line in $\aaa$\dfdot
Next, we will turn the main properties of this operation into axioms.
This leads to the notion of a \emph{ternary ring}.
Our main novelty here is the notation $(a,x,b)\mapsto\langle ax+b \rangle$ for the ternary operation in ternary rings,
replacing the notation $(a,x,b)\mapsto x\cdot a \circ b$ used by M. Hall and other authors.  
The notation $\langle ax+b \rangle$ seems to be much more suggestive than $x\cdot a \circ b$
and makes the whole theory much more transparent.\dfdot

\myitpar{The simplest form of cartesian coordinates on $\aaa$\dfdot}\emph{Fix 
two non-parallel lines $l$\dfcom \dnsp$m$ in $\aaa$\dfdot}
 
The lines $l$\dfcom $m$ allow us to identify $\aaa$ with the cartesian product $l\times m$ in the same way as one does this 
while introducing the cartesian coordinates in the usual Euclidean plane. 
Namely, given a point $p\in\aaa$ we assign to it the pair $(x,y)\in l\times m$\dfcom 
where $x$ is the point of intersection with the line $l$ of the line containing $p$ and parallel to $m$ 
(such a line cannot be parallel to $l$\dfcom because otherwise $l$ would be parallel to $m$ contrary to the assumption), 
and where $y$ is defined in a similar manner. 
This leads to a map $\aaa\rightarrow l\times m$\dfdot 
One can also define a map $l\times m\rightarrow\aaa$ by assigning to $(x,y)$ 
the intersection point of the line containing $x$ and parallel to $m$ and the line containing $y$ and parallel to $l$\dfdot 
Clearly, these two maps are the inverses of each other.
We will use them to identify $\aaa$ with $l\times m$\dfdot
Such an identification is the simplest form of the \emph{cartesian coordinates} on $\aaa$\dfdot

Let $\bf 0$ be the point of intersection of $l$ and $m$\nsp;\sss the point $\bf 0$ serves as the origin of our coordinate system.

\myitpar{Identifying $l$ and $m$\dfdot} Next, we would like to identify the lines $l$ and $m$\dfdot 
In order to do this, we need a line $d$ passing through $\bf 0$ and different from $l$, $m$\dfdot 
One can get such a line, for example, by taking any line $d_0$ intersecting both $l$ and $m$ 
in such a way that the points of intersection are different from $\bf 0$\dfcom
and then take as $d$ the line parallel to $d_0$ and passing through $\bf 0$\dfdot 

\emph{From now on, we will assume that such a line $d$ is fixed.}

Using $d$\dfcom we can construct a natural bijection between $l$ and $m$\dfdot 
Namely, given $x\in l$\dfcom let $z(x)\in d$ be the intersection point with $d$ of the line containing $x$ and parallel to $m$\dfcom 
and let $y(x)\in m$ be the intersection point with $m$ of the line containing $z(x)$ and parallel to $l$\dfdot 
Clearly, $x\mapsto y(x)$ is a bijection $l\rightarrow m$\dfdot

Now, let $K$ be any set endowed with a bijection $K\rightarrow l$\dfdot 
By composing it with the bijection from the previous paragraph, we get a bijection $K\rightarrow m$\dfdot 
Formally, we can simply set $K=l$\dfcom
but we are going to treat $K$ and $l$ differently, and by this reason it is better to consider them as different objects. 
The set $K$ is going to play a role similar to the role of $\bf R$ in Euclidean geometry.

\myitpar{The cartesian coordinates on $\aaa$\dfdot} Our bijections allow us to identify $\aaa$ with $K^2$\dfdot 
We consider this identification as the \emph{cartesian coordinates} on $\aaa$\dfdot
This identification, obviously, turns the line $d$ into the \emph{diagonal}\dss $\{(x,x) : x\in K\}$\dfdot 
Guided by the construction of the usual cartesian coordinates, we denote
the element of $K$ corresponding to the point ${\bf 0}\in l$ by $0$\dfdot 
Then ${\bf 0}=(0,0)$\dfdot

We would like also to have an analogue of the number $1$\dfdot 
In fact, we can choose as $1$ an arbitrary element of $K$ different from $0$\dfdot 
This freedom of choice of $1$ corresponds to the freedom of choice of the unit of measurement in the Euclidean geometry.

\emph{From now on, we will assume that such an element $1\in K$ is fixed.}

\myitpar{The slopes of lines.} Next, we define the \emph{slope} of a line $L$ in $\aaa$\dfdot 
If $L$ is parallel to $l$\dfcom its slope is defined to be $0$\dfdot 
Such lines are called \emph{horizontal}. 
If $L$ is parallel to $m$\dfcom its \emph{slope} is defined to be $\infty$\dfdot 
Such lines are called \emph{vertical}.
If the line $L$ is not vertical, consider the line $L'$ parallel to $L$ and passing through $(0,0)$\dfdot
Let $(1,a)$ be the intersection point of $L'$ with the line $\{(1,z) : z\in K\}$ (
i.e. with the the vertical line passing through $(1,0)$\nsp). 
The \emph{slope} of $L$ is defined to be $a$\dfdot 
Note that it depends on the choice of $1$\dfdot

By the definition, the parallel lines have the same slope.  
Since $d$ contains the point $(1,1)$\dfcom the slope of $d$ is equal to $1$\dfdot 
Clearly, two lines have the same slope if and only if they are parallel.

\myitpar{The ternary operation on $K$\dfdot} Let us define a ternary operation 
$(a,x,b)\longmapsto\langle ax+b \rangle$ on $K$ as follows.
Let $L$ be the unique line intersecting the line $m$ at the point $(0,b)$ and having the slope $a\neq\infty$\dfdot
Since $L$ is a non-vertical line, it intersects the line $m$ at a single point, say $(0,b)$\dfdot  
For every $x\in K$\dfcom the line $L$ intersects the vertical line $\{(x,z) : z\in K\}$ at a unique point.
Let $(x,y)$ be this point and set 
\begin{equation}
\label{line}
\langle ax+b \rangle\qff =\qff y\endss.
\end{equation}
We consider $(a,x,b)\longmapsto\langle ax+b \rangle$ as a ternary operation in $K$\dfdot 
In general we do not have separate multiplication and addition operations in $K$\nsp;\dss 
the angle brackets are intended to stress this.

Clearly, every non-vertical line is the set of points $(x,y)\in K^2$ satisfying (\ref{line}) for some $a,b\in K$\dfdot 
Every vertical line is the set of point $(x,y)\in K^2$ satisfying, for some fixed $c\in K$\dfcom the equation $x=c$\dfdot

\myitpar{The main properties of\trs $(a,x,b)\longmapsto\langle ax+b \rangle$\dfdot} They are the following.
\begin{description}
	\item[T1.] $\langle 1x+0 \rangle\qff =\qff \langle x1+0 \rangle\qff =\qff x$\dfdot
	\item[T2.] $\langle a0+b \rangle\qff =\qff \langle 0a+b \rangle\qff =\qff b$\dfdot
	\item[T3.] \emph{\hspace*{0.4em}If\sss $a,x,y\in K$\dfcom 
	then there is a unique $b\in K$ such that  $\langle ax+b \rangle\qff =\qff y$\dfdot}
	\item[T4.] \emph{\hspace*{0.4em}If\sss $a,a',b,b'\in K$ and $a\neq a'$\dfcom 
	then the equation $\langle ax+b \rangle\qff =\qff \langle a'x+b' \rangle$ has a unique solution $x\in K$\dfdot}
	\item[T5.] \emph{\hspace*{0.4em}If\sss $x,y,x',y'\in K$ and $x\neq x'$\dfcom 
	then there is a unique pair $a,b\in K$ such that $y\qff =\qff \langle ax+b \rangle$ 
	and $y'\qff =\qff \langle ax'+b \rangle$\dfdot}
\end{description}
Let us explain the geometric meaning of these properties. 
This explanation also proves that they indeed hold for $(a,x,b)\longmapsto\langle ax+b \rangle$\dfdot
\begin{description}
  \item[T1\fff:] The equation $\langle 1x+0 \rangle\qff =\qff x$ means that 
  $d\qff =\qff \{(x,x) : x\in K\}$ is a line with the slope $1$\dfdot 
             The equation $\langle x1+0 \rangle\qff =\qff x$ 
             means that the slope of the line passing through $(0,0)$ and $(1,x)$ is equal to $x$ (which is true by the
             definition of the slope). 
  \item[T2\fff:] The equation $\langle a0+b \rangle\qff =\qff b$ 
  means that the line defined by the equation (\ref{line}) intersects $m$ at $(0,b)$
             (which is true by the definition of $\langle ax+b \rangle$). 
             The equation $\langle 0a+b \rangle\qff =\qff b$ means that
             the horizontal line passing through $(0,b)$ consists of points $(a,b)$, $a\in K$\dfdot
  \item[T3\fff:] This means that for every slope $\neq  \infty$ 
  there is a unique line with this slope passing through $(x,y)$\dfdot
  \item[T4\fff:] This means that two lines with different slopes $\neq\infty$ intersect at a unique point.
  \item[T5\fff:] This means that every two points not on the same vertical line (i.e. not on the same line with slope $\infty$) are contained in a unique line with slope $\neq\infty$\dfdot
\end{description}
\myitpar{Ternary rings.} Motivated by these properties, suppose that we have a set $K$ with two distinguished elements $0$ and $1\neq 0$, 
and a ternary operation $(a,x,b)\longmapsto \langle ax+b \rangle$
satisfying {\bf T1}\hff--\hff{\bf T5}. 
Such a $K$ is called a {\em ternary ring}.
Consider the set of points $\aaa=K^2$, and introduce the lines in the following manner: for every $x_0\in K$ we have a line $\{(x_0,y) : y\in K\}$
(such lines are called {\em vertical\/}), 
and for every $a,b\in K$ we have a line $\{(x,y) : y=\langle ax+b \rangle\}$\dfdot 
Clearly, this defines a structure of an affine plane on $\aaa=K^2$ 
(notice that the instances of the axioms of an affine plane involving vertical lines hold trivially).

\mypar{Proposition.}{prop1}  {\em If $K$ is a finite set, then the condition\dss {\bf T5\trs} follows from\dss {\bf T3\trs} and\dss {\bf T4\trs}\nsp\dnsp.}

\proof Given $x,x'\in K$ such that $x\neq x'$\dfcom consider the map $f\colon K^2\rightarrow K^2$ defined by
\[
f(a,b)\qff =\qff (\langle ax+b \rangle, \langle ax'+b \rangle).
\]
Suppose that $f$ is not injective, i.e. that
\begin{equation}
\label{eq1}
\langle ax+b \rangle\qff  =\qff  \langle a'x+b' \rangle, 
\end{equation}
\begin{equation}
\label{eq2}
\langle ax'+b \rangle\qff =\qff \langle a'x'+b' \rangle
\end{equation}
for some $(a,b)\neq (a',b')$\dfdot 
If $a=a'$, then (\ref{eq1}) contradicts {\bf T3}. 
If $a\neq a'$, then the two equalities  (\ref{eq1}), (\ref{eq2}) together contradict {\bf T4}. 
Therefore {\bf T3} and {\bf T4} imply that $f$ is injective. 
Since $f$ is a self-map of a finite set to itself, the injectivity of $f$ implies its surjectivity. 
So, $f$ is a bijection. 
{\bf T5} follows.  \eproof

\mysection{Isomorphisms of ternary rings}{isomorphisms}

\myitpar{The choices involved in the construction of a ternary ring by an affine plane.}
The construction of the ternary ring $K$ associated to an affine plane $\aaa$ involves several choices. 
First, we selected two non-parallel lines $l$ and $m$\dfdot
Then we choose a set $K$ together with a bijection $K\tto l$\dfcom which we may consider as an identification.
This choice is not essential at all.
As we noted, we could simply take $K=l$ and the map $K\tto l$ be the identity map.
 
Then we chose a third line $d$ passing through the intersection point $\bf 0$ of $l$ and $m$\dfcom and a element $1\in K$\dfcom $1\neq 0$\dfdot
The line $d$ and the element $1\in K$ define a point $z\in\aaa$\nsp:\sss 
the intersection point with $d$ of the line parallel to $m$ and passing through the point of $l$ corresponding to $1$\dfcom 
This point corresponds to $(1,1)$ under our identification of $\aaa$ with $K^2$\dfdot 
Conversely, given a point $z\in\aaa$ not on $l,m$\dfcom 
we can define $d$ as the line connecting $\bf 0$ with $z$\dfcom 
and define $1$ as the element of $K$ corresponding to the intersection point with $l$ of the line parallel to $m$ and containing $z$. 
Therefore, the choice of $d$ and $1$ is equivalent to the choice of a point $z$ not contained in the union $l\cup m$\dfdot

\myitpar{The effect of choices.} Let $\aaa '$ be another affine plane with two lines and a point $l',m',z'$ as above,
and let $K'$ be its coordinate ring.
Clearly, there is an isomorphism $f\colon \aaa \rightarrow \aaa '$ 
such that $f(l)=l'$\dfcom $f(m)=m'$\dfcom $f(z)=z'$ if and only if
there is bijection $F\colon K\rightarrow K'$ such that $F(0)=0$\dfcom $F(1)=1$\dfcom and 
\[
F(\langle ax+b \rangle)\qff =\qff \langle F(a)F(x)+F(b) \rangle
\]
for all $a,x,b\in K$\dfdot 
Such a bijection is called an {\em isomorphism} $K\rightarrow K'$\dfdot

\mypar{Proposition.}{prop2} {\em The following two conditions are equivalent.

{\em (i)}\hspace*{0.5em} There is an isomorphism $K^2\rightarrow (K')^2$ taking $\bf 0$ to $\bf 0$\dfcom $K\times 0$ to $K'\times 0$\dfcom
$0\times K$ to $0\times K'$\dfcom and $(1,1)$ to $(1,1)$\dnsp.

{\em (ii)}\hspace*{0.5em} There is an isomorphism $K\rightarrow K'$\dnsp.}

\proof It is sufficient to apply the above observation to $\aaa=K^2$\dfcom $l=K\times 0$\dfcom $m=0\times K$\dfcom $z=(1,1)$\dfcom
and $\aaa=(K')^2$\dfcom $l=K'\times 0$\dfcom $m'=0\times K'$\dfcom $z'=(1,1)$\dfdot  \eproof

We see that up to an isomorphism $K$ is determined by the plane $\aaa$ with a fixed choice of $l,m,z$. We call the
ternary ring $K$ a {\em coordinate ring} of the plane $\aaa$ with a triple $(l,m,z)$ as above.

\mysection{Isotopisms of ternary rings}{isotopisms}

\emph{The later sections do not depend on this one.}

\myitpar{Isotopism.} Ternary rings corresponding to the same affine plane $\aaa$ and the same choice of lines $l\dff,\dff m$\dfcom 
but to different choices of the point $z$\dfcom may lead to non-isomorphic ternary rings.
Still, different choices of $z$ lead to ternary rings which are \emph{isotopic} in the following sense. 
A triple $(F,G,H)$ of bijections $K\rightarrow K'$ is called
an \emph{isotopism}, if $H(0)=0$ and
\[
H(\langle ax+b \rangle)\qff =\qff \langle F(a)G(x)+H(b) \rangle
\]
for all $a,x,b\in K$\dfdot 
Such a triple induces a map $\varphi\colon K^2\rightarrow (K')^2$ by the rule $\varphi (x,y)=(G(x),H(y))$\dfdot 
Clearly, $\varphi$ takes vertical lines to vertical lines. 
The equation $y=\langle ax+b \rangle$ implies $H(y)=\langle F(a)G(x)+H(b) \rangle$\dfcom 
which means that $(x',y')=\varphi (x,y)$ satisfies the equation $y'=\langle F(a)x'+H(b) \rangle$\dfdot 
It follows that $\varphi$ takes the lines with slope $a\neq\infty$ to the lines with slope $F(a)\neq\infty$\dfdot 
We see that $\varphi\colon K^2\rightarrow (K')^2$ is an isomorphism of affine planes.

\mypar{Lemma.}{lemma1} {\em $\varphi$ takes horizontal lines to horizontal lines \textup{(}\fff i.e. $F(0)=0$\dnsp\textup{)},\dss 
and also takes $\bf 0$ to $\bf 0$\dnsp.}

\proof In order to prove the first statement, note that $\varphi$ takes the line $\{(x,0): x\in K\}$ 
to the line $\{(G(x),H(0)): x\in K\}=\{(x',0): x'\in K'\}$
(since $H(0)=0$ and $G$ is a bijection). 
Since both these lines have slope $0$\dfcom we have $F(0)=0$\dfdot 
Therefore, $\varphi$ takes horizontal lines to horizontal lines.
Let $a$, $a'$ be two different slopes. 
Then the lines with equations $y=\langle ax+0 \rangle$\dfcom $y=\langle ax+0 \rangle$ intersect at ${\bf 0}=(0,0)$\dfdot
Their images have the equations $y=\langle F(a)x+0 \rangle$\dfcom $y=\langle F(a')x+0 \rangle$ (recall that $H(0)=0$). 
Since $F$ is a bijection, $F(a)\neq F(a')$ and therefore these two lines intersect only at $\bf 0$\dfdot 
It follows that $\varphi({\bf 0})={\bf 0}$\dfdot  \eproof

\mypar{Corollary.}{cor1} {\em For an isotopism $(F,G,H)$ we have $F(0)=0$ and $G(0)=0$\dfcom in addition to $H(0)=0$\dnsp.}

\proof $F(0)=0$ is already proved. $G(0)=0$ follows from the following two facts: (i) $\varphi({\bf 0})={\bf 0}$\nsp;\dss (ii)
$\varphi$ takes the vertical line $x=0$ to the vertical line $x=G(0)$.  \eproof

\mypar{Corollary.}{cor2} {\em The isomorphism of affine planes induced by an isotopism of ternary rings 
takes the horizontal (respectively, vertical) line containing $\bf 0$ 
to the horizontal (respectively, vertical) line containing $\bf 0$\dnsp.}  \eproof

\mypar{Theorem.}{theorem1} {\em Let $K$ be a coordinate ring of the plane $\aaa$ with a choice of $l\dff,\dff m\dff,\dff z$ as above, and
let $K'$ be the coordinate ring of the plane $\aaa'$ with a choice of $l'\dff,\dff m'\dff,\dff z'$\dfdot 
There is an isomorphism $\aaa\rightarrow\aaa'$ taking $l$ to $l'$ and $m$ to $m'$ \textup{(}\fff but not necessarily $z$ to $z'$\nsp\textup{)} if
and only if there is an isotopism $K\rightarrow K'$\dnsp.}

\proof The {\em ``if''} direction is already proved. Let us prove the {\em ``only if''} direction.

Let us identify $\aaa$ with $K^2$ and $\aaa'$ with $(K')^2$\dfdot
Let $G\colon K\rightarrow K'$ be the map corresponding to the map $K\times 0\rightarrow K'\times 0$ induced by $\varphi$\dfdot 
Similarly, let $H\colon K\rightarrow K'$ be the map corresponding to the map $0\times K\rightarrow 0\times K'$ induced by $\varphi$\dfdot 
Using the fact that every point is the intersection of a unique vertical line with a unique horizontal line, 
and the fact that $\varphi$ maps the vertical (respectively, horizontal) lines to the vertical (respectively, horizontal) lines, 
we see that $\varphi$ is determined by the maps $G\dff,\dff H$, and, in fact, $\varphi (x,y)=(G(x), H(y))$\dfdot 
Clearly, $G(0)=0$ and $H(0)=0$\dfdot

In order to define $F$\dfcom consider for each $a\in K$ the line in $\aaa$ with the slope $a$ passing through $\bf 0$\dfdot 
The map $\varphi$ takes it to a line in $\aaa'$ passing through $\bf 0$\dfdot
Let $F(a)\in K'$ be its slope.

Let us check that $(F,G,H)$ is an isotopism. 
Since $\varphi$ takes parallel lines to parallel lines, 
$\varphi$ takes any line with the slope $a$ to a line with the slope $F(a)$\dfdot 
So, it takes the line with the equation $y=\langle ax+b \rangle$ into a line with the equation of the form $y'=\langle F(a)x'+b' \rangle$\dfdot
The first line contains the point $(0,b)$ (since $\langle a0+b \rangle=b$ by {\bf T2}). 
Therefore, the second line contains the point $\varphi (0,b)=(0,H(b))$\dfdot 
This implies that $H(b)=\langle F(a)0+b'\rangle$\dfdot
But $\langle F(a)0+b'\rangle=b'$ by {\bf T2}. 
Therefore, $b'=H(b)$\dfdot

We see that $\varphi$ maps the line with the equation 
$y=\langle ax+b \rangle$ into the line with the equation $y'=\langle F(a)x'+H(b) \rangle$\dfdot
Since $\varphi (x,y)=(G(x),H(y))$\dfcom we see that $y=\langle ax+b \rangle$ implies $H(y)=\langle F(a)G(x)+H(b) \rangle$\dfdot 
It follows that $(F,G,H)$ is an isotopism.  \eproof

\myitpar{Remark.} If $(F,G,H)$ is an isotopism, then $F$ and $G$ are determined by $H$ and two elements $F^{-1}(1)$\dfcom $G^{-1}(1)$\dfdot
Indeed,
\[
F(a)\qff 
=\qff \langle F(a)1+0 \rangle\qff 
=\qff  \langle F(a)G(G^{-1}(1))+H(0) \rangle\qff 
=\qff H(\langle aG^{-1}(1)+0 \rangle)\endss,
\]
and
\[
G(a)\qff 
=\qff \langle 1G(a)+0 \rangle\qff 
=\qff  \langle F(F^{-1}(1))G(a)+H(0) \rangle\qff 
=\qff H(\langle F^{-1}(1)a+0 \rangle)\endss.
\]

\myitpar{Historical note.} For non-associative algebras, the notion of an equivalence weaker 
than an isomorphism was first introduced by A. A.  Albert \cite{Al}. 
He called two algebras $A$\dfcom $A'$ {\em isotopic} if there is a triple of linear maps $P\dff,\dff Q\dff,\dff R\colon A\tto A'$
such that
\[
R(xy)\qff =\qff P(x)Q(y)\endss.
\]
He called such a triple an {\em isotopy} of $A$ and $A'$\dfdot 
Albert relates that
\begin{quote}
\emph{The concept of isotopy was suggested to the author by the work of N. Steenrod who,
in his study of homotopy groups in topology, was led to study isotopy of division algebras.}
\end{quote}
Albert noticed that if associativity of the multiplication is not assumed, the notion of isotopy
is more suitable than the obvious notion of isomorphism, which leads to too many non-isomorphic 
(but isotopic) algebras.

It is only natural that the notion of an isomorphism is too narrow for the ternary rings also. 
The corresponding notion of an {\em isotopism} was introduced by
M. V. D. Burmester \cite{Bu}, and, independently, by D. Knuth \cite{Kn}. 
Both Burmester and Knuth proved Theorem \ref{theorem1} above. 
D. Knuth \cite{Kn}, moreover, found an affine plane $\aaa$ such that all ternary rings corresponding
to different choices of $z$ (but the same choice of $l\dff,\dff m$\nsp) are pairwise non-isomorphic. 
His plane is finite, and the corresponding ternary rings have $32$ elements. 
See \cite{Kn}, Section 5. 
Unfortunately, his plane was found with the help of a computer, and, as Knuth writes,
 \emph{``No way to construct this plane, except by trial and error, is known.''}
To the best knowledge of the author, this is still the case.

\mysection{Veblen-Wedderburn systems}{quasi-fields}

\myitpar{The left Veblen-Wedderburn systems.} Let $K$ be a set with two binary operations $(x\fff,\fff y)\dnsp\mapsto x+y$ and $(x\fff,\fff y)\mapsto x\fff y$\dnsp,\dss 
called the \emph{addition} and the \emph{multiplication}, respectively,
 and two distinguished elements $0$\dfcom $1$\dfcom $0\neq 1$\dfdot
If the following properties\trs {\bf VW1}\fff--\fff{\bf VW5}\trs hold, 
$K$ is called a\sss \emph{left Veblen-Wedderburn system},\sss 
or, more recently, a\sss \emph{left quasi-field}. 
\begin{description}
  \item[VW1.] \emph{$K$ is an abelian group with respect to the addition $+$\ffdot}
  \item[VW2.] \emph{\hspace*{0.4em}Given $a\dff,\dff b\neq 0$\dfcom each of the equations $a\dff x=b$ and $x\dff a=b$ has a unique solution $x$\nsp;\sss \hspace*{1.0em}moreover, this solution is $\neq 0$\dfdot In addition, if $a,b\neq 0$\dfcom then $ab\neq 0$\dfdot}
  \item[VW3.] \emph{$1\dff x=x\dff 1=x$\dfcom $0\dff x=x\dff 0=0$\dfcom and $x+0=0+x=x$ for all $x$\dfdot}
  \item[VW4.] \emph{\hspace*{0.4em}Left distributivity: $a(x+y)=ax+ay$ for all $a\dff,\dff x\dff,\dff y$\dfdot}
  \item[VW5.] \emph{\hspace*{0.4em}For $a\neq a'$\dfcom the equation $a\dff x=a'\hff x+b$ has a unique solution $x$\dfdot}
\end{description}
This notion was introduced by O. Veblen and J. Wedderburn \cite{VW}. 
Notice that\trs {\bf VW5}\trs is a weak version of the right distributivity. 
Clearly, under conditions\trs {\bf VW1},\trs {\bf VW2}\trs it follows from the right distributivity.

\myitpar{The right Veblen-Wedderburn systems.} In order to define \emph{right Veblen-Wedderburn system}, or \emph{right quasi-fields}, 
we replace\trs {\bf VW4}\trs and\trs {\bf VW5}\trs by the following two conditions. 
\begin{description}
  \item[VW4-r{}.] \emph{Right distributivity: $(x+y)a=x\fff a+y\fff a$ for all $a\dff,\dff x\dff,\dff y$\dfdot}
  \item[VW5-r{}.] \emph{For $a\neq a'$\dfcom the equation $x\fff a=x\fff a'+b$ has a unique solution $x$\dfdot}
\end{description}
Clearly, $K$ is a right quasi-field if and only if $K$ with the same addition, $0$\dfcom $1$\dfcom and the opposite multiplication
$a\cdot b=b\fff a$\dfcom is a left quasi-field.

\myitpar{Weak quasi-fields.} If $K$ satisfies only conditions\trs {\bf VW1}\fff--\fff{\bf VW4},\trs 
it is called a \emph{a weak left quasi-field}. 
Similarly, $K$ is called a \emph{weak right quasi-field}, if it satisfies conditions\trs {\bf VW1}\fff--\fff{\bf VW3}\trs and\trs {\bf VW4-r{}.}

\myitpar{From Veblen-Wedderburn systems to ternary rings.} If $K$ is a left or right quasi-field, 
then we can define a ternary operation $(a\dff,\dff x\dff,\dff b)\mapsto \langle a\dff x+b \rangle$ by the obvious rule $\langle a\fff x+b \rangle=a\fff x+b$\dfdot 
We claim that $K$ with this ternary operation and the distinguished elements $0$ and $1$ is a ternary ring.
Let us check this first for left quasi-fields.
\begin{description}
  \item[T1\fff:] This condition follows from\trs {\bf VW3}.
  \item[T2\fff:] This condition also follows from\trs {\bf VW3}.
  \item[T3\fff:] This condition follows from\trs {\bf VW1}.
  \item[T4\fff:] Let $a\dff,\dff a'\dff,\dff b\dff,\dff b'\in K$ and $a\neq a'$\dfdot 
  The equation $a\fff x+b =a'\hff x+b'$ for $x$ is equivalent to $a\fff x=a'\hff x+(b'-b)$ by\trs {\bf VW1}. 
  It has a unique solution by\trs {\bf VW5}.
  \item[T5\fff:] Let $x\dff,\dff y\dff,\dff x'\dff,\dff y'\in K$ and $x\neq x'$\dfdot 
  The equations $y=a\fff x+b$ and $y'=a\fff x'+b$ for $a\dff,\dff b$ imply
  \[
  y-y'\qff =\qff a\fff x-a\fff x'
  \]
  by\trs {\bf VW1}, and hence imply
  \[
  y-y'\qff =\qff a(x-x')
  \]
  by\trs {\bf VW4}. 
  If $y\neq y'$\dfcom this equation is uniquely solvable for $a$ by\trs {\bf VW2}. 
  If we know $a$\dfcom we can find $b$ from
  either of the equations $y=a\fff x+b$\dfcom $y'=a\fff x'+b$\dfdot 
  Therefore $b$ is unique. 
  This proves\trs {\bf T5} in the case $y\neq y'$\dfdot 
  If $y=y'$\dfcom then $a$ has to be equal to $0$ by\trs {\bf VW2}\trs (\fff since $x-x'\neq 0$\nsp). 
  Therefore $b=y=y'$\dfdot 
  This proves\trs {\bf T5} in the case $y=y'$\dfdot 
\end{description}
For a right quasi-field $K$ the conditions\trs {\bf T1}\fff--\fff{\bf T3}\trs hold by the same reasons as for the left quasi-fields 
(they do not depend on the distributivity). 
Let us check\trs {\bf T4}\trs and\trs {\bf T5}.
\begin{description}
  \item[T4\fff:] Let $a\dff,\dff a'\dff,\dff b\dff,\dff b'\in K$ and $a\neq a'$\dfdot 
   The equation $ax+b =a'x+b'$  for $x$ is equivalent to $(a-a')x=(b'-b)$ by\trs {\bf VW1} and\trs 
  {\bf VW4-r{}} (the right distributivity). 
  It has a unique solution by\trs {\bf VW2} and\trs {\bf VW3}\dss (the latter is needed if $b'-b=0$).
  \item[T5\fff:] Let $x\dff,\dff y\dff,\dff x'\dff,\dff y'\in K$ and $x\neq x'$\dfdot 
  The equations $y=ax+b$ and $y'=ax'+b$ for $a\dff,\dff b$ imply
  \[
  ax\qff =\qff ax'+(y-y')\endss.
  \]
  Since $x\neq x'$\dfcom this equation has a unique solution $a$ by {\bf VW5-r}. 
  As above, if we know $a$\dfcom we can find $b$ from either of the equations $y=ax+b$\dfcom $y'=ax'+b$\dfdot 
  Therefore $b$ is unique. 
  This proves \trs{\bf T5}.
\end{description}
Notice that going from left to right quasi-fields switches the roles of\trs {\bf VW4}\trs and\trs {\bf VW5}.

\myitpar{Reconstructing quasi-field from the corresponding ternary ring.} A left quasi-field can be restored 
from the corresponding ternary ring in an obvious manner: it has the same $0$ and $1$\nsp;
the addition and the multiplication are defined by $a+b=\langle 1a+b \rangle$ and $ab=\langle ab+0\rangle$\dfdot 
Indeed, $1(ax+0)+b=ax+0+b=ax+b$\dfdot 
Therefore, we may consider quasi-fields as a special class of ternary rings. 
In particular, a quasi-field $K$ defines an affine plane. 
Of course, this plane can be described directly: its set of points is $K^2$\dfcom 
and its lines are given by the equations of the form $x=a$ 
and of the form $y=ax+b$\dfcom where $(x,y)\in K^2$ and $a\dff,\dff b$ are fixed elements of $K$\dfdot

\mypar{Proposition.}{prop3}  \emph{If\trs $K$ is weak left quasi-field and is finite, 
then $K$ is a left quasi-field \textup{(}\fff i.e.\trs {\bf VW5}\trs follows
from\trs {\bf VW1}\fff--\fff{\bf VW4\trs}\trs if\trs $K$\dss is finite\fff\textup{)}.}

\proof For $a\neq a'$\dfcom let $f(x)=ax-a'x$\dfdot 
Suppose that $f$ is not injective, i.e. $ax-a'x=ay-a'y$ for some $x\neq y$\dfdot 
Then $a(x-y)=a'(x-y)$ by\trs {\bf VW1}\trs and\trs {\bf VW4}\trs (the left distributivity). 
Since $a\neq a'$\dfcom this contradicts\trs {\bf VW2}. 
Therefore, $f$ is injective. 
Being a self-map of a finite set to itself, it is bijective (cf. the proof of Proposition \ref{prop1}). 
Therefore, for every $b$ there is a unique $x$ such that  $ax-a'x=b$\dfdot 
Hence,\trs {\bf VW5}\trs holds.  \eproof

\myitpar{Finiteness.} Proposition \ref{prop1} shows that in the finite case we can drop\trs {\bf T5}\trs form the axioms of a ternary ring. 
By Proposition \ref{prop3} we can also drop\trs {\bf VW5}\trs form the axioms of a quasi-field for finite $K$\dfdot 
While checking\trs {\bf T4\trs} for the ternary ring associated to a quasi-field above, we referred to\trs {\bf VW5}. 
If the quasi-field is finite and we drop the axiom\trs {\bf VW5}, we have to use the Proposition \ref{prop3}, 
and the role of\trs {\bf VW5}\trs is passed to the left distributivity.\\

In some situations the finiteness can be replaced by the finite dimensionality over an appropriate skew-field.

\mypar{Proposition.}{prop4} \emph{Suppose that a weak left quasi-field $K$ contains a subset $F$ 
which is a skew-field with respect to the same operations and with the same $0$ and $1$\dfdot 
Suppose that, in addition, 
\begin{equation}
\label{weak-1}
(x\dff y)a\qff =\qff x\fff(y\dff a),
\end{equation}
\begin{equation}
\label{weak-2}
(x+y)a\qff =\qff x\dff a+y\dff a\endss,
\end{equation}
for all $a\in F$ and $x\dff,\dff y\in K$\dfdot 
Then $K$ is a right vector space over $F$\dfdot 
If this vector space is finitely dimensional, then $K$ is a left quasi-field \textup{(}i.e. the condition\trs {\bf VW5}\trs holds\textup{)}.}

\proof The first statement is clear. 

Let us prove the second one. 
For $a\in K$\dfcom let $L_a\colon K \rightarrow K$ be the left multiplication by $a$\dfcom 
i.e. $L_a(x)=ax$\dfdot 
By\trs {\bf VW4}\trs we have $L_a(x+y)=L_a(x)+L_a(y)$ for all $x\dff,\dff y\in K$\dfdot 
Moreover, if $b\in F$\dfcom then $L_a(x\dff b)=a(x\dff b)=(a\dff x)b=L_a(x)b$\dfdot
It follows that $L_a$ is (right) linear map of the vector space $K$ to itself.

We need to check that for $a\neq a'$ the equation $L_a(x)=L_{a'}(x)+b$ has a unique solution $x$\dfdot
Let $L=L_a-L_{a'}$\ffdot
It is sufficient to show that the equation $L(x)=b$ has a unique solution $x$\dfdot 
Clearly, $L$ is a linear map. 
If $L(y)=0$\dfcom then $a\dff y-a'\hff y=0$ and $a\dff y=a'\hff y$\dfdot 
Since $a\neq a'$\dfcom the condition\trs {\bf VW2}\trs implies that this is possible only if $y=0$\dfdot 
We see that $L$ is linear self-map of $K$ with trivial kernel. 
Since $K$ is assumed to be finitely dimensional, $L$ is an isomorphism.
This implies that $L(x)=b$ has a unique solution. 
This proves the second statement of the proposition.  \eproof

Our proof of Proposition \ref{prop4} follows the proof of Theorem 7.3 in \cite{HP}.

\mysection{Near-fields, skew-fields, and isomorphisms}{near-fields}

In general, if affine planes $K^2$ and $(K')^2$ are isomorphic, the ternary rings $K$ and $K'$ do not need to be isomorphic. 
The goal of this section is to prove that they will be isomorphic if $K'$ is a skew-field. 
See Corollary \ref{cor3} below. 
A part of the proof works in a greater generality, namely for near-fields, which we will define in a moment.

A \emph{left near-field} is a left quasi-field with associative multiplication. 
Non-zero elements of a left near-field form a group with respect to the multiplication. 
The \emph{right near-fields} are defined in an obvious manner. 
Clearly, being a skew-field is equivalent to be being a left and right near-field simultaneously.

\mypar{Lemma.}{lemma2} \emph{Let\trs $K'$\dss be a left near-field. 
Let\dss ${\bf 0}=(0,0)\in (K')^2$\dfcom and let\trs $l\dff,\dff m$\dss be, respectively, 
the horizontal and the vertical lines in\dss $(K')^2$\sss passing through\dss $\bf 0$ \textup{(}i.e. $l=K'\times 0$\sss and\dss $m=0\times K'$\nsp\textup{)}\hff. 
For every two points\sss $z\dff,\dff z'\in (K')^2$\sss not in\dss $l\cup m$\dfcom
there is an automorphism of the affine plane $(K')^2$\sss preserving\sss $\bf 0$\dfcom $l$\dfcom and\sss $m$\dfcom and taking $z$\sss to\sss $z'$\dfdot}

\proof It is sufficient to consider the case when $z=(1,1)$\dfdot 
Let $z'=(u\dff,\dff v)$\dfdot 
Since $(u\dff,\dff v)$ is not on $l\dff,\dff m$\dfcom both $u$ and $v$ are non-zero. 
Consider the map $f\colon (K')^2\rightarrow (K')^2$ defined by $f(x\dff,\dff y)=(ux\dff,\dff vy)$\dfdot 
Clearly, $f(1\dff,\dff 1)=(u\dff,\dff v)$\dfcom and $f$ takes the vertical line $x=a$ to the vertical line $x=au$\dfdot 
If $y=ax+b$\dfcom then $vy=v(ax)+vb=(vau^{-1})ux+vb$ by the left distributivity and the associativity of the multiplication 
(here, as usual, $u^{-1}$ is the unique solution of the equation $xu=1$\dnsp)\hff. 
It follows that $f$ takes the line $y=ax+b$ to the line $y=(vau^{-1})x+vb$\dfdot 
Hence $f$ is an automorphism of $(K')^2$\dfdot  \eproof

\mypar{Corollary.}{cor3} {\em Let $K$ be a ternary ring, and let $K'$ be a left near-field. 
Suppose that there is an isomorphism of planes $f\colon K^2\rightarrow (K')^2$ 
taking $\bf 0$ to $\bf 0$ and taking the horizontal \textup{(}\fff respectively, vertical\fff\textup{)} line through $\bf 0$ in $K^2$
to horizontal \textup{(}\fff respectively, 
vertical\fff\textup{)} line through $\bf 0$ in $(K')^2$\dfdot 
Then $K$ is isomorphic to $K'$ as a ternary ring.}

\proof Let $z=(1\dff,\dff 1)\in K^2$. By taking the composition of $f$ with an appropriate automorphism $g\colon (K')^2\rightarrow (K')^2$\dfcom
if necessary, we can assume that $f(1\dff,\dff 1)=(1\dff,\dff 1)$ (the required $g$ exists by the lemma). 
It remains to apply Proposition \ref{prop2}.  \eproof

\mypar{Lemma.}{lemma3} {\em Let $K'$ be a skew-field.\trs 
Let ${\bf 0}=(0\dff,\dff 0)\in (K')^2$\dfcom and let\dss $l\dff,\dff m$\dss be, respectively, 
the horizontal and the vertical lines in\sss $(K')^2$\sss passing through $\bf 0$ 
\textup{(}\fff i.e. $l=K'\times 0$ and $m=0\times K'$\nsp\textup{)}. 
Let\dss $l'\dff,\dff m'$\dss be any two non-parallel lines in\sss $(K')^2$\dfdot 
Then there is an automorphism of the affine plane\sss $(K')^2$\sss taking\sss $l$\sss to\sss $l'$\sss and\sss $m$\sss to\sss $m'$\sss 
\textup{(}\dff and, in particular, taking\sss $\bf 0$\sss to the intersection point of\dss $l'$\dss and\dss $m'$\nsp\textup{)}\hff.}

\proof If $K$ is a field, this is a fact well-known from the linear algebra. 
In general, one needs to check that there is no need to use commutativity of the multiplication. 
Let us first check that some natural maps are isomorphisms.
\begin{description}
  \item[\rm (i)] The map $D(x\dff,\dff y)=(y\dff,\dff x)$ is an isomorphism. 
  Indeed, it takes the line $x=a$ to the line $y=0x+a$\dfcom and the line $y=0x+b$ to the line $x=b$\dfdot
  If $a\neq 0$\dfcom it takes the line $y=ax+b$\dfcom i.e.
  the line $x=a^{-1}y-a^{-1}b$ (where $a^{-1}$ is the unique solution of the equation $xa=1$\nsp) to the line $y=a^{-1}x-a^{-1}b$\dfdot 
  Here we used the left distributivity and the associativity of the multiplication.
  \item[\rm (ii)] For any $c\dff,\dff d\in K'$ the map $f(x\dff,\dff y)=(x+c\dff,\dff y+d)$ is an isomorphism. 
  Indeed, it takes the line $x=a$ to the line $x=a+c$, and the line $y=ax+b$, to the line $y=ax-ac+b+d$\dfdot
  Here we used the left distributivity.
  \item[\rm (iii)] For any $c\in K'$ the map $f(x\dff,\dff y)=(x\dff,\dff y-cx)$ is an isomorphism. 
  Indeed, it takes every line $x=a$ to itself, and it takes 
  the line $y=ax+b$ to the line $y=(a-c)x+b$\dfdot
  Here we used the right distributivity.
  \item[\rm (iv)] For any $c\in K'$ the map $g(x\dff,\dff y)=(x-cy\dff,\dff y)$ is an isomorphism. 
  Indeed,  $g=D\circ f\circ D$\dfcom where $D(x\dff,\dff y)=(y\dff,\dff x)$ 
  and $f(x\dff,\dff y)=(x\dff,\dff y-cx)$\dfdot
\end{description}
By using an isomorphism of type\dss (ii)\dss if necessary, we can assume that $l'\dff,\dff m'$ intersect at $\bf 0$\dfdot
By using the isomorphism $D$ from\dss (i)\dss if necessary, we can assume that $l'$ is not equal to $m=0\times K$\dfdot 
Then $l'$ has the form $y=cx$\dfdot 
The map $f(x\dff,\dff y)=(x\dff,\dff y-cx)$ is of type\dss (iii)\dss and takes $l'$ to $l$\dfdot 
Therefore we can assume that $l'=l$ and $m'$ intersects $l'=l$ at $\bf 0$\dfdot 
Then $m'$ has a equation of the form $x=cy$ and an isomorphism of type\dss (iv)\dss takes $m'$ to $m$\dfdot 
Since any automorphisms of type\dss (iv)\dss takes $l$ to $l$\dfcom this completes the proof.  \eproof

\mypar{Theorem.}{theorem2} {\em Let $K$ be a ternary ring, and let $K'$ be a skew-field. 
Suppose that there is an isomorphism of planes $f\colon K^2\rightarrow (K')^2$\dfdot 
Then $K$ is isomorphic to $K'$ as a ternary ring.}

\proof This follows from Lemma \ref{lemma3} and Corollary \ref{cor3}.  \eproof

It follows that in order to construct an affine plane not coming from a skew-field, it is sufficient to construct
a quasi-field which is not a skew-field 
(a quasi-field which is isomorphic to a skew-field is a skew-field itself). 
In particular, it is sufficient to construct a left quasi-field in which 
the right distributivity does not hold or a right quasi-field in which the left distributivity does not hold.
Alternatively, it is sufficient to construct a (left or right) quasi-field in which
the associativity of multiplication does not hold. 
We will present a construction of such quasi-field in Section \ref{andre}.

\mysection{Translations}{translations}

The previous section provided us with a method of constructing affine planes not isomorphic to any affine plane defined by a skew-field.
In this section we will present another method, based on an investigation of special automorphisms of affine planes called \emph{translations}.
It allows to show that some planes are not isomorphic even to any plane defined by a left quasi-field (see Theorem \ref{theorem3}\trs below).

Let $\aaa$ be an affine plane. An automorphism $f\colon\aaa\rightarrow\aaa$ is called a \emph{translation} 
if $f(l)$ is parallel to $l$ for every line $l$ (equal lines are considered to be parallel), 
and if $f$ preserves every line from a class of parallel lines. 
Clearly, for a non-trivial (i.e., not equal to the identity) translation there is exactly one such class of parallel lines. 
Every line from this class is called a\sss \emph{trace}\sss of\sss $f$\dfdot 
If $\aaa$ is realized as $K^2$ for a ternary ring $K$\dfcom 
then a translation is called\dss \emph{horizontal}\dss if it preserves all horizontal lines, 
i.e. if the class of horizontal lines is its trace. 

\myitpar{\hspace*{-0.3em}}
\emph{The next two propositions are not used in the rest of the paper.}

\mypar{Proposition.}{prop5} \emph{A non-trivial translation has no fixed points.}

\proof Let $f$ be a translation fixing a point $z$\dfdot 
Let $l$ be trace of $f$\dfdot 

Let $m$ be a line containing $z$ and not parallel to $l$\dfdot 
Since $f(m)$ is parallel to $m$ and contains $z$\dfcom we have $f(m)=m$\dfdot 
Every point of $m$ is the unique intersection point of $m$ and a line parallel to $l$\dfdot 
The map $f$ leaves both of these lines invariant. 
Therefore $f$ fixes all points of $m$\dfdot 

We see that $f$ fixes all points except, possibly, the points of the line $l_z$ passing through $z$ and parallel to $l$\dfdot
By applying the same argument to any point not on $l_z$ in the role of $z$\dfcom 
we conclude that $f$ fixes the points of $l_z$ also, i.e. that $f={\rm id}$\dfdot  \eproof

\mypar{Proposition.}{prop6} \emph{Let $z,z'$ be two different points, and let $l$ be the line passing through $z\dff,\dff z'$\dfdot 
There is no more than one translation taking $z$ to $z'$\dfcom 
and if such a translation exists, it leaves invariant every line parallel to $l$\dfdot}

\proof If $f_1\dff,\dff f_2$ are two different translations taking $z$ to $z'$\dfcom 
then $f_1^{-1}\circ f_2$ is a non-trivial translation fixing $z$\dfcom 
contradicting to Proposition \ref{prop5}. 

Now, let $f$ be a translation such that $f(z)=z'$\dfcom 
and let $m$ be a trace of $f$\dfdot 
Let $m_z$ be the line passing through $z$ and parallel to $m$\dfdot 
Then $m_z$ is also a trace of $f$\dfdot 
Clearly, we have $z\in m_z$ and $z'=f(z)\in m_z$\dfdot 
It follows that $l=m_z$ and, hence, $l$ is a trace of $f$\dfdot 
Therefore, $f$ leaves invariant every line parallel to $l$\dfdot 
This completes the proof.  \eproof

\mypar{Lemma.}{lemma4} \emph{Let $K$ be a left quasi-field. 
For every two points $(c_1\dff,\dff d_1)$\dfcom $(c_2\dff,\dff d_2)$ of the plane $K^2$\dfcom 
there is a translation of $K^2$ taking $(c_1\dff,\dff d_1)$ to\sss $(c_2\dff,\dff d_2)$\dfdot}

\proof In this case there are obvious maps expected to be translations, 
namely the maps of the form $f(x\dff,\dff y)=(x+c\dff,\dff y+d)$\dfcom where $c\dff,\dff d\in K$\dfdot 
Clearly, if $c=c_2-c_1$\dfcom $d=d_2-d_1$\dfcom then $f(c_1\dff,\dff d_1)=(c_2\dff,\dff d_2)$\dfdot 
Let us check that these maps are indeed translations.

The map $f(x\dff,\dff y)=(x+c\dff,\dff y+d)$ takes the line $x=a$ to the line $x=a+c$\dfcom 
and the line $y=ax+b$ to the line $y=ax-ac+b+d$\dfdot
(Cf. (ii) in the proof of Lemma \ref{lemma3} in the previous section.) 
In particular, it takes vertical lines to vertical lines,
and the lines with the slope $a$ to the lines with the slope $a$\dfdot 
If $c=0$\dfcom then $f$ preserves all vertical lines, and therefore is a translation. 
If $c\neq 0$\dfcom then $d=ec$ for some $e$ by {\bf VW2}. 
Since $f$ takes the line $y=ex+b$ to the line $y=ex-ec+b+d$ and $ex-ec+b+d=ex-d+b+d=ex+b$\dfcom 
we see that $f$ leaves invariant every line with the slope $e$\dfdot 
It follows that $f$ is a translation in this case also.  \eproof

\mypar{Lemma.}{lemma5} \emph{Let $K$ be a right quasi-field. 
Suppose that for every $v\in K$\dfcom $v\neq 0$\dfcom 
the plane $K^2$ admits a translation taking $\bf 0$ to $(v,0)$\dfdot 
Then the left distributivity law holds in $K$\dnsp.}

\proof Let $f$ be a translation such that $f(0\dff,\dff 0)=(v\dff,\dff 0)$\dfdot 
Since the line passing through $(0\dff,\dff 0)$ and $(v\dff,\dff0)$ is the horizontal line $K\times 0$\dfcom 
the translation $f$ is a horizontal translation. 
Let us show that $f$ has the expected form $f(a\dff,\dff b)=(a+v\dff,\dff b)$\dfdot
This follows from the following four observations.
\newpage
\begin{description}
  \item[\rm 1.] The line $y=x$ with the slope $1$ passing through $(0\dff,\dff 0)$ is mapped to the line with the slope $1$ passing through 
  $(v\dff,\dff 0)$\dfcom i.e. to the line $y=x-v$\dfdot
  \item[\rm 2.] For every $a\in K$, the map $f$ preserves the line $y=a$\dfdot 
  It follows that $f$ takes the point of intersection
  of the lines $y=a$ and $y=x$ to the point of intersection of the lines $y=a$ and $y=x-v$\dfdot 
  This means that $f(a\dff,\dff a)=(a+v\dff,\dff a)$\dfdot
  \item[\rm 3.] The vertical line $x=a$ containing $(a\dff,\dff a)$ is mapped to the vertical line containing $(a+v\dff,\dff a)$\dfcom 
  i.e. to the line $x=a+v$\dfdot
  \item[\rm 4.] The point of intersection of the lines $y=b$ and $x=a$ is mapped to the point of intersection of the lines
  $y=b$ and $x=a+v$\dfdot 
  In other terms, $f(a\dff,\dff b)=(a+v\dff,\dff b)$\dfdot
\end{description}
Now, $f$ takes the line $y=ax$ containing $(0\dff,\dff 0)$ to another line with the slope $a$\dfcom 
i.e. to a line of the form $y=ax-c$\dfdot 
Since it contains $(v,0)$, we have $c=av$\dfdot 
So, the line $y=ax$ is mapped to the line $y=ax-av$\dfdot 
For every $u\in K$ the point $(u,au)$ belongs to the line $y=ax$ and is mapped to the point $(u+v\dff,\dff au)$\dfdot 
Therefore, $(u+v\dff,\dff au)$ belongs to the line $y=ax-av$\dfcom i.e. $au=a(u+v)-av$\dfcom or $a(u+v)=au+av$\dfdot

Since this is true for all $a\dff,\dff u\dff,\dff v\in K$\dfcom the left distributivity holds.  \eproof

\mypar{Theorem.}{theorem3} \emph{Let $K$ be a right quasi-field for which the left distributivity does not hold, 
and let $K'$ be a a left quasi-field. 
Then the planes $K^2$ and $(K')^2$ are not isomorphic.}

\proof By Lemma \ref{lemma5} there is a point $(c\dff,\dff d)\in K^2$ 
such that no translation takes $(0\dff,\dff 0)$ to $(c,d)$\dfdot 
Therefore,  $K^2$ is not isomorphic to any plane constructed from a left quasi-field by Lemma \ref{lemma4}.  \eproof

\mysection{Andr\'e quasi-fields}{andre}

\myitpar{The norm map.} Let $K$ be a field, and let $G$ be a finite group of automorphisms of $K$\dfdot 
Let $F$ be the subfield of $K$ consisting of all elements fixed by $G$\dfdot
By the Galois theory, the dimension of $K$ as a vector space over $F$ is equal to the order of $G$\nsp; 
in particular, it is finite. 
The \emph{norm map} $N$ is defined as follows:
\[
N(x)\qff =\qff \prod_{g\in G} g(x)\endss.
\]
Clearly, $N(x)\in F$ for all $x\in K$\dfdot 
Moreover, $N$ defines a homomorphism $K^*\rightarrow F^*$ from the multiplicative group $K^*$ to the multiplicative group $F^*$ of the
fields $K$\dfcom $F$ respectively. 
Obviously, $N(g(a))=N(a)$ for any $g\in G$\dfdot

\myitpar{Modifying the multiplication.} Note that $N(1)=1$\dfcom and, therefore, $1\in N(K^*)$\dfdot
Let $\varphi\colon N(K^*)\rightarrow G$ be a map subject to the only condition $\varphi(1)=1$\dfdot 
In particular, $\varphi$ does not need to be a homomorphism. 
Given such a map $\varphi$\dfcom we construct a new multiplication $\odot$ in $K$ as follows. 
Of course, the new multiplication $\odot$ will depend on $\varphi$\dfcom but we will omit this dependence from the notations.
Let $\alpha$ be equal to $\varphi\circ N$ on $K^*$\dfcom and let $\alpha(0)=1={\rm id}_K$\dfdot 
So, $\alpha$ is a map $K\rightarrow G$\dfdot 
We will often denote $\alpha(x)$ by $\alpha_x$\nsp; it is an automorphism $K\rightarrow K$ belonging to the group $G$\dfdot
The multiplication $\odot$ is defined by the formula
\[
x\odot y\qff  =\qff x\dff\alpha_x(y)\endss,
\]
for all $x\dff,\dff y\in K$\dfdot
Let $K_{\varphi}$ be the set $K$ endowed with the same addition and the same elements $0$\dfcom $1$ as $K$\dfcom and with the multiplication $\odot$\dfdot

\mypar{Theorem.}{theorem4} \emph{$K_{\varphi}$ is a left quasi-field.}

\proof {\bf VW1\fff:} This property holds for $K_{\varphi}$ because it holds for $K$\dfdot

{\bf VW3\fff:} Note that $\alpha_1=\varphi (N(1))=\varphi(1)=1$\dfdot 
This implies $1\odot x=x$ for all $x$\dfdot 
Also, since $\alpha_x$ is an automorphism of $K$\dfcom
we have $\alpha_x(1)=1$\dfcom $\alpha_x(0)=0$\dfcom and, therefore, $x\odot 1=x$\dfcom $x\odot 0=0$ for all $x\in K$\dfdot 
In addition, $0\odot y=0\dff\alpha_0(y)=0\dff y=0$\dfdot 
These observations imply the multiplicative part of {\bf VW3} for $K_{\varphi}$\nsp;\trs 
the additive part holds for $K_{\varphi}$ because it holds for $K$\dfdot 

{\bf VW4\fff:} Since $\alpha_x$ is an automorphism of $K$\dfcom 
we have $\alpha_x(y+z)=\alpha_x(y)+\alpha_x(z)$\dfcom 
and therefore $x\odot (y+z)=x\odot y + x\odot z$\dfdot
So, the left distributivity law {\bf VW4} holds for $K_{\varphi}$\dfdot

{\bf VW2\fff:} Suppose that $a\dff,\dff b\neq 0$\dfdot 
First of all, notice that $\alpha_a(b)\neq 0$ 
(because $\alpha_a$ is an automorphism of $K$\nsp), and, therefore, $a\odot b \neq 0$\dfdot 
Next, consider the equation $a\odot x=b$\dfdot 
It is equivalent to $a\dff\alpha_a(x)=b$\dfcom which, in turn,
is equivalent to $\alpha_a(x)=a^{{\minus}1}b$\dfdot 
It follows that $x=\alpha_a^{{\minus}1}(a^{{\minus}1}b)$ is the unique solution. 

It remains to consider the equation $x\odot a=b$\dfdot 
Notice that 
\[
N(x\odot a)\qff =\qff N(x\fff\alpha_x(a))\qff =\qff N(x)N(\alpha_x(a))\qff =N\qff (x)\dff N(a)\endss,
\] 
since $N(g(a))\qff  =\qff  N(a)$ for any $g\in G$\dfdot
Therefore,\qss $x\odot a\qff  =\qff  b$ implies that $N(x)\dff N(a)\qff  =\qff  N(b)$\dfdot 
This, in turn, implies that 
\[
N(x)\qff =\qff N(a)^{{\minus}1}N(b)\qff =\qff N(a^{{\minus}1}b)\endss,
\] 
and 
$\displaystyle
\alpha(x)\qff =\qff \varphi(N(x))\qff =\qff \varphi(N(a^{-1}b))\qff =\qff \alpha(a^{{\minus}1}b)$.
In other terms, $\alpha_x\qff =\qff \alpha_{a^{{\minus}1}b}$\dfdot 
\newpage 

Therefore,\qss if\qss $x\dff\alpha_x(a)\qff =\qff x\odot a\qff =\qff b$\dfcom 
then $x\dff\alpha_{a^{{\minus}1}b}(a)\qff =\qff b$ and 
\begin{equation}
\label{solution}
x\qff =\qff b(\alpha_{a^{{\minus}1}b}(a))^{{\minus}1}.
\end{equation}
It follows that the equation $x\odot a\qff =\qff b$ has no more than one solution. 

Let us check that (\ref{solution}) is, indeed, a solution.
Let $g=\alpha_{a^{-1}b}$\dfdot 
If $x$ is defined by (\ref{solution}), then
\[
\alpha_x\qff =\qff \alpha(x)\qff =\qff \varphi(N(b(\alpha_{a^{-1}b}(a))^{-1})\qff 
=\qff \varphi(N(b\fff (g(a)^{{\minus}1}))\endss.
\]
At the same time,
\begin{align*}
  N(b\fff (g(a)^{{\minus}1}))  
&= N(b)\dff N(g(a)^{{\minus}1})\hspace*{1em}\mbox{ (\fff because $N$ is a homomorphism) } \\
&= N(b)\dff N(g(a^{{\minus}1}))\hspace*{1em}\mbox{ (\fff because $g\in G$ and hence $g(a)^{-1}=g(a^{-1})$\dff) } \\
&= N(b)\dff N(a^{{\minus}1})\hspace*{2.4em}\mbox{ (\fff because $g\in G$ and hence $N(g(b))=N(b)$\dff) } \\
&= N(ba^{{\minus}1})\dff\hspace*{4.0em}\mbox{ (\fff because $N$ is a homomorphism) } \\
&= N(a^{{\minus}1}b)\dff\hspace*{4.0em}\mbox{ (\fff because the multiplication in $K$ is commutative)\endss. }
\end{align*}
It follows that
$\displaystyle
\alpha_x=\varphi(N(b\fff (g(a)^{{\minus}1}))=\varphi(N(a^{{\minus}1}b))=\alpha_{a^{{\minus}1}b}\endss.
$
Therefore
\[
x\odot a=x\dff\alpha_x(a)=x\dff\alpha_{a^{{\minus}1}b}(a)=b\dff (\alpha_{a^{{\minus}1}b}(a))^{{\minus}1}\dff\alpha_{a^{{\minus}1}b}(a)=b
\]
This proves that (\ref{solution}) is a solution of $x\odot a=b$ and  completes our verification of {\bf VW2}.

{\bf VW5\fff:} It remains to check {\bf VW5}. 
To this end, we will apply Proposition \ref{prop4} (see Section \ref{quasi-fields}).
We already established that $K_{\varphi}$ is a weak left quasi-field. 
Notice that 
\begin{equation}
\label{alphaxa}
\alpha_x(a)=a\hspace*{1.9em}\mbox{\fff for\hspace*{0.5em} any\hspace*{0.5em} $x\in K$\hspace*{0.5em} and\hspace*{0.5em} $a\in F$}\endss,
\end{equation} 
because $F$ is fixed by all elements of $G$\dfdot 
Therefore, 
\begin{equation}
\label{odot}
x\odot a\qff  =\qff  x\dff a\hspace*{1.5em}
\mbox{\fff for\hspace*{0.5em} all\hspace*{0.5em} $x\in K$,\hspace*{0.5em} $a\in F$}\endss.
\end{equation} 
This immediately implies the condition (\ref{weak-2}) of Proposition \ref{prop4}. 
Now, let $x,y\in K$\dfcom and $a\in F$\dfdot 
The following calculation shows that the condition (\ref{weak-1}) of Proposition \ref{prop4} also holds:
\begin{align*}
(x\odot y) \odot a\qff  & =\qff (x\odot y)\dff a\hspace*{2.3em}\mbox{ (\fff by (\ref{odot})\dff) } \\
& =\qff x\dff\alpha_x(y)\dff a\hspace*{2.3em}\mbox{ (\fff by the definition of $\odot$\fff) } \\
& =\qff x\dff\alpha_x(y)\fff\alpha_x(a)\hspace*{0.5em}\mbox{ (\fff by (\ref{alphaxa})\dff) }  \\
& =\qff x\dff\alpha_x(y a)\hspace*{2.3em}\mbox{ (\fff because $\alpha_x$ is a homomorphism) } \\
& =\qff x\odot (y a) \hspace*{2.3em}\mbox{ (\fff by the definition of $\odot$\fff) } \\
& =\qff x\odot (y\odot a)\hspace*{1.15em}\mbox{ (\fff by (\ref{odot})\dff). }
\end{align*}
\newpage

Since $K$ is finitely dimensional vector space over $F$\dfcom the Proposition \ref{prop4} applies. 
It implies that\dss {\bf VW5}\trs holds and hence $K_{\varphi}$ is a left quasi-field. 
This completes the proof.  \eproof

The left quasi-fields $K_{\varphi}$ are called {\em left Andr\'e quasi-fields}. 
The {\em right Andr\'e quasi-fields} are
constructed in a similar manner, with the multiplication given by the formula $x\odot y = \alpha_y(x)\dff y$\dfdot 
As the next two theorems show, in a left Andr\'e quasi-field
the multiplication is almost never associative (Theorem \ref{theorem5}), and the right distributivity holds only if $\varphi$ is the
\emph{trivial map}, i.e. the map taking every element to $1\in G$ (Theorem \ref{theorem6}). 
Of course, the corresponding results hold for the right Andr\'e quasi-fields.
We will call an Andr\'e quasi-field {\em non-trivial} if $\varphi$ is a non-trivial map.

\mypar{Theorem.}{theorem5} {\em The multiplication in\dss $K_{\varphi}$ is associative 
if and and only if $\varphi$ is a homomorphism $N(K^*)\rightarrow G$\dfcom 
i.e. if and only if $\varphi(u\dff v)=\varphi(u)\dff\varphi(v)$ for all $u\dff,\dff v\in N(K^*)$\dfdot}

\proof Let $x\dff,\dff y\in K^*$ and $g\in G$\dfdot
Then
\begin{align*}
\alpha(x\dff g(y))\qff  & =\qff \varphi(N(x\dff g(y))\hspace*{4.5em}\mbox{\fff (by the definition of $\alpha$) } \\
 & =\qff \varphi(N(x)\dff N(g(y)))\hspace*{2.5em}\mbox{\fff (because $N$ is a homomorphism) } \\
 & =\qff \varphi(N(x)\dff N(y))\hspace*{3.9em}\mbox{\fff (because $N(y)\qff =\qff N(g(y))$) } \\
 & =\qff \varphi(N(x\dff y))\qff =\qff \alpha(x\dff y).
\end{align*}
Therefore, $\alpha(x\dff g(y))\qff  =\qff  \alpha(x\dff y)$\dfcom or, equivalently, 
$\alpha_{x\dff g(y)}\qff  =\qff  \alpha_{x\dff y}$ for all $x\dff,\dff y\in K$ and $g\in G$\dfdot
{\qff}It follows that
\begin{equation}
\label{xyz}
\alpha_{x\dff g(y)}\dff (z)\qff  =\qff  \alpha_{x\dff y}\dff (z)
\end{equation}
for all $x\dff,\dff y\in K^*$\dfcom $z\in K$\dfcom and $g\in G$\dfdot

If $g=\alpha_x$\dfcom then $x\dff g(y)\qff =\qff x\dff\alpha_x(y)\qff =\qff x\odot y$ 
and hence (\ref{xyz}) turns into
\begin{equation}
\label{axy}
\alpha_{x\odot y}\dff (z)\qff  =\qff  \alpha_{x\dff y}\dff (z)\endss.
\end{equation}
By applying (\ref{axy}) we can compute $(x\odot y)\odot z$ as follows: 
\[
(x\odot y)\odot z\qff  
=\qff (x\odot y)\dff\alpha_{x\odot y}(z)\qff  
=\qff (x\odot y)\dff\alpha_{x\dff y}(z)\qff 
=\qff  x\dff\alpha_x(y)\dff\alpha_{x\dff y}(z).
\]
Next, let us compute $x\odot (y\odot z)$:  
\[
x\odot (y\odot z)\qff  
=\qff x\odot (y\dff\alpha_y(z)) 
=\qff x\dff\alpha_x(y\dff\alpha_y(z)) 
=\qff x\dff\alpha_x(y)\dff\alpha_x(\alpha_y(z)).
\]
By comparing the results of these computations, we see that for $x\dff,\dff y\in K^*$ 
the associativity law $(x\odot y)\odot z = x\odot (y\odot z)$ holds if and only if
$\alpha_{x\dff y}(z) = \alpha_x(\alpha_y(z))$\dfdot
Since the associativity law trivially holds when $x=0$ or $y=0$\dfcom
the associativity law for $\odot$ holds if and only if $\alpha_{x\dff y}(z) = \alpha_x(\alpha_y(z))$
for all $x\dff,\dff y\dff,\dff z\in K$\dfcom
or, equivalently, if and only if $\alpha_{x\dff y} =\alpha_x\circ\alpha_y$ for all $x\dff,\dff y\in K$\dfdot

Recalling the definition of $\alpha$\dfcom we see that $\alpha_{x\dff y}\qff  =\qff \alpha_x\circ\alpha_y$ is equivalent to
\[
\varphi(N(x))\dff\varphi(N(y))\qff =\qff \varphi(N(x\dff y))\endss,
\]
and since $N(x\dff y)\qff =\qff N(x)\dff N(y)$\dfcom is equivalent to
\begin{equation}
\label{phi}
\varphi(N(x))\dff\varphi(N(y))\qff =\qff \varphi(N(x)\dff N(y))\endss.
\end{equation}
Clearly, (\ref{phi}) holds for all $x\dff,\dff y\in K$ if and only if $\varphi\colon N(K^*)\tto G^*$ is homomorphism.
It follows that the associativity is equivalent to $\varphi$ being a homomorphism.  \eproof

\mypar{Theorem.}{theorem6} {\em The right distributivity law holds for $K_{\varphi}$ if and only if $\varphi$ maps all
elements of $N(K^*)$ to $1\in G$ \textup{(}\fff and therefore $K_{\varphi}=K$\nsp\textup{)}.}

\proof Clearly, if $\varphi$ maps $N(K^*)$ to $1$, then $K_{\varphi}=K$ and the right distributivity holds.

Suppose now that the right distributivity holds. Then for all $x\dff ,\dff y\dff ,\dff a\in K$ we have
\[
x\odot a + y\odot a\qff  =\qff  (x+y)\odot a\endss.
\]
We can rewrite this as
\begin{equation}
\label{dist-1}
x\dff X(a)+y\dff Y(a)\qff =\qff (x+y)\dff Z(a)\endss,
\end{equation}
where $X=\alpha_x$\dfcom $Y=\alpha_y$\dfcom $Z=\alpha_{x+y}$\dfdot 
By using the fact that $X\dff ,\dff Y\dff ,\dff Z$ are automorphisms of $K$
and applying (\ref{dist-1}) to $a\hff b$ in the role of $a$\dfcom we get
\begin{equation}
\label{dist-2}
x\dff X(a)\dff X(b)+y\dff Y(a)\hff Y(b)\qff =\qff x\dff X(a\hff b)+y\dff Y(a\hff b)
\end{equation}
\[
\phantom{x\dff X(a)\dff X(b)+y\dff Y(a)\hff Y(b)\qff }
=\qff (x+y)\dff Z(a\hff b)\qff =\qff (x+y)\dff Z(a)\dff Z(b)\endss.
\]
By combining (\ref{dist-2}) with (\ref{dist-1}), we get
\[
x\dff X(a)\dff X(b)+y\hff Y(a)\dff Y(b)\qff =\qff (x\dff X(a)+y\hff Y(a))\dff Z(b)\endss.
\]
Let us multiply this identity by $x+y$, and then apply (\ref{dist-1}) to $b$ in the role of $a$\nsp:
\begin{equation*}
  (x+y)(x\dff X(a)\dff X(b)+y\hff Y(a)\hff Y(b))\qff
  \end{equation*} 
  \begin{equation*} 
  =\qff (x+y)(x\dff X(a)+y\hff Y(a))\dff Z(b)\qff 
  \end{equation*} 
  
  \vspace*{-2.7\bigskipamount}
  \begin{equation*}
   \hspace*{0em} 
  =\qff (x+y)\dff Z(b)\dff (x\dff X(a)+y\hff Y(a))\qff 
  =\qff (x\dff X(b)+y\hff Y(b))(x\dff X(a)+y\hff Y(a))\endss.
\end{equation*}

By opening the parentheses and canceling the equal terms, we get
\[
y\dff x\dff X(a)\dff X(b)+x\dff y\hff Y(a)\hff Y(b)\qff =\qff x\dff y\dff X(b)\hff Y(a)+y\dff x\hff Y(b)\dff X(a)\endss.
\]
Suppose that $x\dff ,\dff y\neq 0$\dfdot 
Then we can divide the last equation by $x\dff y\neq 0$ and get
\[
Y(a)\hff Y(b)+X(a)\dff X(b)\qff =\qff X(b)\hff Y(a)+Y(b)\dff X(a)\endss.
\]
The last identity is equivalent to
\[
Y(a)\hff Y(b)-Y(a)\dff X(b)+X(a)\dff X(b)-X(a)\hff Y(b)\qff =\qff 0\endss,
\]
and, therefore, to
\begin{equation}
\label{yxab}
(Y(a)-X(a))\hff (Y(b)-X(b))\qff =\qff 0\endss.
\end{equation}
This holds for all $a\dff ,\dff b\in K$\dfdot 
If $Y(a)-X(a)\neq 0$ for some $a$\dfcom 
then (\ref{yxab}) implies that $Y(b)-X(b)=0$ for all $b$\dfcom
and, in particular, for $b=a$ in contradiction with  $Y(a)-X(a)\neq 0$\dfdot
It follows that $X(a)=Y(a)$ for all $a$\dfdot 
In other terms, $X=Y$\dfdot 
By recalling that $X=\alpha_x$\dfcom $Y=\alpha_y$\dfcom 
and that $x\dff ,\dff y$ are arbitrary non-zero elements of $K$\dfcom 
we conclude that all automorphisms $\alpha_x$ with $x\neq 0$ are equal,
and, in particular, are equal to $\alpha_1$\dfdot 
But the $\alpha_1=\varphi(N(1))=\varphi(1)=1$ by the assumption.
It follows that $\varphi(N(x))=\alpha_x=1$ for all $x\in K^*$\dfcom 
and hence $\varphi(z)=1$ for all $z\in N(K^*)$\dfdot 
This completes the proof.  \eproof

The Galois theory provides many explicit examples of field $K$ with a finite group of automorphisms $G$. The freedom
of choice of the map $\varphi$ allows to construct left Andr\'e quasi-field with non-associative multiplication
(by using Theorem \ref{theorem5}), and left Andr\'e quasi-field in which the right distributivity does not hold (by using Theorem \ref{theorem6}).
One can also construct a left Andr\'e quasi-field with associative multiplication in which the right distributivity does not hold.
We leave this as an exercise for the readers moderately familiar with Galois theory.

\myitpar{Remark.} In this section we to a big extent followed the exposition in \cite{HP}, Section IX.3.

\mysection{Conclusion: non-Desarguesian planes}{conclusion}

If $K$ is a left quasi-field with non-associative multiplication (say, a left Andr\'e quasi-field), then $K$ is not isomorphic to any skew-field. 
By Theorem \ref{theorem2}, $K^2$ is not isomorphic to any plane defined over a skew-field, and, 
therefore, is a non-Desarguesian plane.\newpage

If $K$ is a left quasi-field in which the right distributivity does not hold, then, again, $K$ is not isomorphic to any skew-field. 
By Theorem \ref{theorem2}, $K^2$ is a non-Desarguesian plane.

Let $K$ be a right quasi-field in which the left distributivity does not hold. 
For example, one can take as $K$ any nontrivial right Andr\'e quasi-field 
(we can take as $K$ a left Andr\'e quasi-field with the opposite multiplication, 
or use the right Andr\'e quasi-field version of Theorem \ref{theorem6}). 
Then, by Theorem \ref{theorem3}, $K^2$ is not isomorphic to any plane defined over a left quasi-field. 
In particular, $K^2$ is not isomorphic to any plane defined over a skew-field, and, therefore, is a non-Desarguesian plane.

\begin{flushright} 

October 8, 2014 

April 15, 2016 (minor edits)

\vspace*{\bigskipamount}

http:/\!/\hspace*{-0.07em}nikolaivivanov.com
\end{flushright}


\begin{thebibliography}{XXX}



\bibitem[Bu]{Bu} M. V. D. Burmester,\qss On the non-unicity of translation planes over division algebras, {\em Archive Math.},
Vol. XV (1964), 364--370.

\bibitem[Al]{Al} A. A. Albert,\qss Non-associative algebras: I. Fundamental concepts and isotopy, {\em Annals of Math.},
Vol. 43, No. 4. (1942), 685--707.

\bibitem[Ar]{Ar} E. Artin,\qss \emph{Geometric algebra}, Wiley, Hoboken, NJ, 1988; Reprint of the 1st edition, Interscience Publishers,  New York, 1957.

\bibitem[H1]{Ha1} M. Hall, Jr.,\qss \emph{Theory of groups}, AMS Chelsea Publishing, 1999; 1st edition: Macmillam, New York, 1959.

\bibitem[H2]{Ha2} M. Hall, Jr.,\qss \emph{Combinatorial theory}, Wiley, Hoboken, NJ, 1986.

\bibitem[Ha]{Har} R. Harstshorne,\qss \emph{Foundations of projective geometry}, W.A. Benjamin, 1967.

\bibitem[HP]{HP} D. R. Hughes, F. C. Piper,\qss \emph{Projective planes}, Graduate Texts in Mathematics, V. 6, Springer-Verlag, 1973.

\bibitem[Kn]{Kn} D. Knuth,\qss Finite semifields and projective
planes, {\em J. Algebra}, V. 2, No. 2 (1965), 182--217.

\bibitem[VW]{VW} O. Veblen and J. Wedderburn,\qss Non-Desarguesian and non-Pascalian geometries,
{\em Trans. AMS}, V. 8 (1907), 379--388.

\bibitem[W]{W} Ch. Weibel,\qss Survey of non-Desarguesian planes, {\em Notices of the AMS}, V. 54. No 10 (2007), 1294--1303.  



\end{thebibliography}
\end{document}